\documentclass[11pt,twoside]{amsart}
\usepackage{amsmath, amsthm, amscd, amsfonts, amssymb, graphicx, color}
\usepackage[bookmarksnumbered, colorlinks, plainpages]{hyperref}
\usepackage{pgf,tikz}
\usetikzlibrary{arrows}
\usepackage[utf8]{inputenc}
\usepackage{pstricks-add}
\input{mathrsfs.sty}
\newtheorem{thm}{Theorem}[section]

\newtheorem{rem}[thm]{\bf{Remark}}
\newtheorem{example}[thm]{\bf{Example}}

\numberwithin{equation}{section}

\begin{document}


\title{On well-f-coveredness of lexicographic products of graphs}
\author{Reza Jafarpour-Golzari}
\address{Department of Mathematics, Payame Noor University, P.O.BOX 19395-3697 Tehran, Iran; Department of Mathematics, Institute for Advanced Studies
in Basic Science (IASBS), P.O.Box 45195-1159, Zanjan, Iran}

\email{r.golzary@iasbs.ac.ir}

\thanks{{\scriptsize
\hskip -0.4 true cm MSC(2010): Primary: 05C70; Secondary: 05C38, 05C76.
\\
\newline\indent{\scriptsize}}}

\maketitle


\begin{abstract}
A simple graph $G$ is said to be well-f-covered, whenever any two maximal induced forest in $G$ be of the same order. In this note, well-f-coveredness of lexicographic product of two graphs in case where the first component is empty, is characterized. In cases where the second component is empty, and the second component is nonempty, a necessary condition is given, and in each one, by an example, it is shown that the given condition is not sufficient.
\end{abstract}

\vskip 0.8 true cm


\section{\bf Introduction}
\vskip 0.4 true cm
In the sequel, we refer to \cite{Wes} for any backgrounds in graph theory. Also all graphs are considered to be finite and simple.

A graph $G$ is called well-f-covered, whenever all its maximal induced forests be of the same order. The concept of well-f-coveredness was introduced by the author for the first time, and some properties of well-f-covered graphs and their behavior under several graph operations was studied (see \cite{Jaf}).

In this paper, the well-f-coveredness of lexicographic product of two graphs in case where the first component is empty, is characterized (Theorems 3.1). In cases where the second component is empty, and the second component is nonempty, a necessary condition is given (Theorems 3.2 and 3.5) and in each one, by an example, it is shown that the given condition is not sufficient (Examples 3.4 and 3.7). The author hopes that the current study to be completed in the future with characterizing the well-f-coveredness of lexicographic product of graphs in general. Also it is hoped that this paper provides motivation to study well-coveredness of other types of product in graphs. 

\vskip 0.8 true cm

\section{\bf Preliminaries}
\vskip 0.4 true cm
Let $G$ be a graph. For a subset $X$ of $V(G)$, we denote the subgraph of $G$ induced by $X$, by $G[X]$. For convenience, in the sequel, we will assume that subgraph means induced subgraph. 

 In the graph $G$, a subset $X$ of $V(G)$ is called an independent vertex set, whenever no two vertices $x$ and $y$ in $X$ be adjacent. The size of a maximum independent set in $G$ is said to be the independence number of $G$ and is denoted by $\alpha (G)$. The graph $G$ is called well-covered, if all its maximal independent sets be of the same size \cite{Plu1}.
 
 A graph $G$ is said to be a forest, whenever $G$ does not contain any cycle. A leaf in a forest is any vertex of degree 1. Every nontrivial forest has at least two leaves. 
 
The forest number of  a graph $G$, denoted by $f(G)$, is the order of a maximum forest in $G$. A graph $G$ is called well-f-covered, whenever all its maximal forests be of the same order $f(G)$ (see \cite{Jaf}). 
 
 \vskip 0.8 true cm

\section{\bf The main results}
\vskip 0.4 true cm
First of all, it can be seen that for every two empty graphs $G$ and $H$, $GoH$ is empty and therefore is well-f-covered.\\

The following theorem characterizes well-f-coveredness of the lexicographic product of two graphs where the first component is empty.

\begin{thm}
Consider an empty graph $G$ with $m$ vertices. For any graph $H$, $GoH$ is well-f-covered, if and only if $H$ be well-f-covered, and anyway $f(GoH)=mf(H)$.
\begin{proof}
We use induction on $m$.

\textbf{Base case:} Suppose that $m=1$. In this case, it is cleat that $GoH\cong H$. Therefore $GoH$  is well-f-covered, if and only if $H$ be well-f-covered, and anyway $f(GoH)$=1f(H).

\textbf{Induction step:} Suppose that the statement holds for $m(\in\mathbb{N})$. let $G$ be an empty graph with $m+1$ vertices and $H$ be any graph. Let the vertices of $G$ be $u_{1}, \ldots, u_{m}, u_{m+1}$ and let $G_{1}$ and $G_{2}$ be the empty graphs with vertex sets $\{u_{1}, \ldots, u_{m}\}$ and $\{u_{m+1}\}$, respectively. Since $G$ is empty, $GoH$ is the union of two disjoint graphs $G_{1}oH$ and $G_{2}oH$. Now, if $GoH$ be well-f-covered, by Theorem 6.1 in \cite{Jaf}, $G_{1}oH$ and $G_{2}oH$ both are well-f-covered and $f(GoH)=f(G_{1}oH)+f(G_{2}oH)$; Therefore $H$ is well-f-covered and
\[f(G_{1}oH)=mf(G), \ \ \ f(G_{2}oH)=1f(H).\]
Thus $f(GoH)(m+1)f(H)$. Conversely, if $H$ be well-f-covered, $G_{1}oH$ and $G_{2}oH$ are well-f-covered and
\[f(G_{1}oH)=mf(G), \ \ \ f(G_{2}oH)=1f(H).\] 
Again by Theorem 6.1 in \cite{Jaf}, $GoH$ is well-f-covered and
\[f(GoH)=f(G_{1}oH)+f(G_{2}oH)=mf(H)+1f(H)=(m+1)f(H).\]
\end{proof}
\end{thm}

On the other hand, if $H$ be a trivial graph, then for any graph $G$, one has $GoH\cong G$, and therefore $GoH$ is well-f-covered, if and only if $G$ be well-f-covered, and anyway $f(GoH)=f(G)$.

The following theorem gives a necessary condition for well-f-coveredness of lexicographic product of two graphs in case where the second component is empty.

\begin{thm}
Let $G$ be a graph and $H$ be an empty graph of order $n$. If $GoH$ be well-f-covered, then for every maximal forest F in $G$,
\[n(I(F)+K_{2}(F)+L(F))+K_{2}(F)+L^{'}(F)=f(GoH)\]
where $I(F)$ is the number of isolated vertices in $F$, $K_{2}(F)$ is the number of connected components of $F$ which are $K_{2}$, $L(F)$ is the number of leaves in the other components of $F$, and $L^{'}(F)$ is the number of vertices of degree at least 2 in $F$.
\begin{proof}
Suppose that $GoH$ be well-f-covered and consider an arbitrary maximal forest $F$ in $G$. Set:
\vskip 0.3 true cm
\hskip -0.3 true cm\small{$X$}\tiny{$:=\{g\in V(G)|\ $$g$ is an isolated vertex of $F$ or a leaf in a component (of $F$) which is not $K_{2}.\},$}
\vskip 0.3 true cm
\hskip 2 true cm\small{$Y$}\tiny{$:=\{g\in V(G)|\ $$g$ is a vertex of degree at least 2 in $F.\}.$} 
\vskip 0.3 true cm
\hskip -0.4 true cm\normalsize{Also,} take a vertex from any component of $F$ which is $K_{2}$ and denote by $Z$, the set of these vertices, and also denote by $T$, the set of other vertices of such components. Set:
\[V^{*}:=((X\cup Z)\times V(H))\cup ((Y\cup T)\times \{h\})\]
where $h$ is a fixed vertex in $H$. We claim that $(GoH)[V^{*}]$ is a maximal forest of $GoH$. Therefore since $GoH$ is well-f-covered,
\[f(GoH)=|V^{*}|=n(|X|+|Z|)+|Y|+|T|\]
\[\hskip 3.5 true cm =n(I(F)+L(F)+K_{2}(F))+L^{'}(F)+K_{2}(F).\]

\textbf{Proof of the claim:} Let $(GoH)[V^{*}]$ contains a cycle as:
\[(g_{1},h_{1})(g_{2},h_{2}), \ldots, (g_{k-1},h_{k-1})(g_{k},h_{k})(=(g_{1},h_{1})).\]
If $g_{1}$ be an isolated vertex in $F$, then $g_{1}=g_{2}$ and $h_{1}$ and $h_{2}$ are adjacent in $H$, a contradiction. If $g_{1}$ be a leaf in a component of $F$ which is not $K_{2}$, then $g_{2}\in Y$ and therefore $h_{2}=h$. Since $g_{1}$ is of degree 1 in $F$, $g_{k-1}=g_{2}$ and therefore $h_{k-1}=h$. Thus the vertex $(g_{2},h_{2})$ is repeated in the cycle, a contradiction. If $g_{1}\in Y$, no one of the vertices $g_{2}, \ldots, g_{k-1},g_{k}$ can not be a leaf. Thus $h_{2}= \cdots= h_{k-1}=h_{k}=h$ and therefore the vertices $g_{i}$, $1\leq i\leq k-1$, are distinct pairwise (note that the vertices $(g_{i},h_{i})$ are distinct pairwise as vertices of $GoH$ in the cycle). Therefore
\[g_{1}g_{2} \ldots g_{k-1}g_{k}(=g_{1})\] 
is a cycle in $F$, a contradiction. Finally, if $g_{1}\in T$, one have $g_{3}=g_{1}$ and therefore $(g_{1},h)$ is repeated, a contradiction again. Thus $(GoH)[V^{*}]$ is a forest.

Let $(a,b)\in V(GoH)$ does not be in the forest $(GoH)[V^{*}]$. If $a\notin V(F)$, by maximality of $F$, the subgraph $G[V(F)\cup \{a\}]$ contains a cycle as:
\[(a=)a_{1}a_{2} \ldots a_{l-1}a_{l}(=a_{1})\]
and therefore
\[(a,b)(a_{2},h), \ldots, (a_{l-1},h)(a_{l},h)(=(a,b))\]
is a cycle in $(GoH)[V^{*}\cup \{(a,b)\}$, a contradiction. If $a\in V(F)$, then $a$ is in $Y\cup T$ and also $b\neq h$. Now, if $a$ be in $Y$, $a$ is adjacent with two distinct vertices $a^{'}$ and $a^{''}$ in $F$ and
\[(a^{''},h)(a,b)(a^{'},h)(a,h)(a^{''},h)\]
is a cycle in $(GoH)[V^{*}\cup \{(a,b)\}]$, and if $a\in T$, $a$ adjacent with a leaf $a^{'''}\in Z$ in $F$ and 
\[(a,b)(a^{'''},h)(a,h)(a^{'''},b)(a,b)\]
is a cycle, a contradiction.
\end{proof}
\end{thm}

\begin{rem}
In Theorem 3.2, if $n=1$, the well-f-coveredness of $GoH$ implies that for every maximal forest $F$ in $G$, 
\[f(GoH)=I(F)+K_{2}(F)+L(F)+K_{2}(F)+L^{'}(F)=|V(F)|.\] 
This means that $G$ is well-f-covered and $f(G)=f(GoH)$, something we are already had as well.
\end{rem}

The following example shows that the condition given in Theorem 3.2, is not a sufficient condition. 

\begin{example}
Consider the graph $P_{4}o2K_{1}$. Every maximal forest $F$ of $P_{4}$ is an induced $P_{3}$ and therefore
\[2(I(F)+K_{2}(F)+L(F))+K_{2}(F)+L^{'}(F)=6=f(P_{4}o2K_{1}).\]
But the composition is not well-f-covered.
\end{example}

Now we give a necessary condition in the case where the second component is not empty.

\begin{thm}
For a nonempty graph $G$ and a nonempty graph $H$, if $GoH$ is well-f-covered, then  \\
(1) $G$ is well-covered, and if every maximal forest of $G$ be without isolated vertex and $H$ has a maximal independent set of size 1, then $G$ is well-f-covered and $f(G)=f(GoH),$ \\
(2) $H$ is well-f-covered, and if $G$ has a maximal forest with at least one leaf, then $H$ is well-covered, \\
(3) $f(GoH)=\alpha (G)f(H)$, \\
(4) If $F$ be a maximal forest in $G$, then for any maximal independent set $M_{H}$ in $H$,
\[f(H)I(F)+|M_{H}|(K_{2}(F)+L(F))+K_{2}(F)+L^{'}(F)=f(GoH),\]
where $I(F)$, $K_{2}(F)$, $L(F)$, and $L^{'}(F)$ are on pair with Theorem 3.2.
\begin{proof}
Let $GoH$ is well-f-covered. We show that (1), (2), (3), and (4) are satisfied. 

Consider a fixed maximal forest $F_{H}$ in $H$. Let $M$ be any maximal independent set in $G$. Set:
\[V_{M}:=M\times V(F_{H}).\]
If $(GoH)[V_{M}]$ has a cycle:
\[(x_{1},y_{1})(x_{2},y_{2}) \ldots (x_{n-1},y_{n-1})(x_{n},y_{n})(=(x_{1},y_{1})),\]
then $x_{1}=x_{2}= \cdots x_{n-1}=x_{n}$ and
\[y_{1}y_{2} \cdots y_{n-1}y_{n}\]
is a cycle in $F$, a contradiction. Thus $(GoH)[V_{M}]$ is a forest. We show hat this forest is maximal in $GoH$. Let $(v,w)\notin V_{M}$ be a vertex of $GoH$. If $v\notin M$, then by maximality of $M$, there exists a vertex of $M$, say $a$, adjacent with $v$ and therefore for two adjacent vertices $w_{1}$ and $w_{2}$ in $F$ (note that such two vertices are exist because $H$ is not empty),
\[(a,w_{1})(v,w)(a,w_{2})(a,w_{1})\]
is a cycle in $(GoH)[V_{M}\cup \{(v,w)\}]$, and if $v\in M$, then $w\notin V(F_{H})$ and therefore adding $w$ to $F_{H}$ forms a cycle in $H$ as:
\[l_{1}l_{2} \cdots l_{k-1}l_{k}(=l_{1})\]
and therefore
\[(v,l_{1})(v,l_{2}) \ldots (v,l_{k-1})(v,l_{k})(=(v,l_{1}))\]
is a cycle in  $(GoH)[V_{M}\cup \{(v,w)\}]$. By well-f-coveredness of $GoH$, we have:
\[|M|=\frac{f(GoH)}{|V(F_{H})|}.\]
Since $\frac{f(GoH)}{|V(F_{H})|}$ is independent of $M$, $G$ is well-covered.

Now, for any two maximal forests $F_{1H}$ and $F_{1H}$ in $H$,
\[|F_{1H}|=|F_{2H}|=\frac{f(GoH)}{\alpha (G)}\]
and therefore $H$ is well-f-covered and
\[f(GoH)=\alpha (G)f(H).\]

Suppose that $F$ be a maximal forest in $G$. Define $X$, $Y$, $Z$, and $T$ the same as Theorem 3.2 and define $X_{1}$ to be the set of all isolated vertices of $F$, and $X_{2}:=X\smallsetminus X_{1}$. Let $M_{H}$ and $F_{H}$, be an independent set and a maximal forest in $H$, respectively. Set:
\[V^{*}:=(X_{1}\times V(F_{H}))\cup ((X_{2}\cup Z)\times M_{H}))\cup ((Y\cup T)\times \{h\})\]
where $h$ is a fixed vertex in $M_{H}$. We claim that $(GoH)[V^{*}]$ is a maximal forest in $GoH$. Therefore since $GoH$ is well-f-covered and $X_{1}$, $X_{2}\cup Z$, and $Y\cup T$ are disjoint pairwise, the equality in (4) holds.

\textbf{Proof of the claim:} Let $(GoH)[V^{*}]$ contains a cycle as:
\[(g_{1},h_{1})(g_{2},h_{2}), \ldots, (g_{p-1},h_{p-1})(g_{p},h_{p})(=(g_{1},h_{1})).\]
If $g_{1}$ be an isolated vertex of $F$, then $g_{1}=g_{2}= \cdots g_{p-1}=g_{p}$, and therefore $h_{1}h_{2} \cdots h_{p-1}h_{p}(=h_{1})$ is a cycle in $F_{H}$, a contradiction. Thus $g_{1}$ is in $X_{2}\cup Y\cup Z\cup T$. Similarly $g_{2}, \ldots, g_{p-1}, g_{p}\in X_{2}\cup Y\cup Z\cup T$. Hence for every $1\leq j\leq p$, $h_{j}\in M_{H}$. If for a $1\leq j\leq p-1$, $g_{j}=g_{j+1}$, then $h_{j}$ and $h_{j+1}$ are adjacent, something which is not possible; Therefore $g_{j}$ and $g_{j+1}$ are adjacent for every  $1\leq j\leq p-1$, and 
\[g_{1}g_{2} \cdots g_{p-1}g_{p}(=g_{1})\]
is a closed path in $F$. Now, in each one of the cases $g_{1}\in X_{2}$, $g_{1}\in Y$, $g_{1}\in Z$, and $g_{1}\in T$, similarly to what we saw in the proof of Theorem 3.2, a contradiction is revealed. Hence $(GoH)[V^{*}]$ is a forest.

Consider a vertex $(\alpha,\beta)$ in $GoH$ such that $(\alpha,\beta)\notin V^{*}$. We show that $(GoH)[V^{*}\cup \{(\alpha,\beta)\}]$ contains a cycle. If $\alpha\notin V(F)$, since $F$ is maximal, there is a cycle as:
\[(\alpha =)u_{1}u_{2} \cdots u_{t-1}u_{t}(=u_{1})\]
in $G[V(F)\cup \{\alpha\}]$; Therefore
\[((\alpha,\beta)=)(u_{1},h)(u_{2},h), \ldots, (u_{t-1},h)(u_{t},h)(=(u_{1},h))\]
is a cycle in $(GoH)[V^{*}\cup \{(\alpha,\beta)\}]$. Let $\alpha\in F$. If $\alpha\in X_{1}$, then $\beta\notin V(F_{H})$; By minimality of $f_{H}$, there is a cycle in $H[V(F)\cup \{\beta\}]$, and by adding the first component $\alpha$ to each one of the vertices of this cycle, a cycle in $(GoH)[V^{*}\cup \{(\alpha,\beta)\}]$ is formed. If $\alpha\in X_{2}$, then $\alpha$ is adjacent with a vertex in $Y$, say $\alpha^{'}$; Since $(\alpha,\beta)\notin V^{*}$, $\beta\notin M_{H}$ and therefore $\beta$ is adjacent with a vertex $\beta^{'}\in M_{H}$; Thus
\[(\alpha,\beta)(\alpha^{'},h)(\alpha,\beta ^{'})(\alpha,\beta)\]
is a cycle. If $\alpha\in Z$, then $\alpha$ is adjacent with a leaf in $T$, say $\alpha^{''}$, and 
\[(\alpha,\beta)(\alpha^{''},h)(\alpha,\beta ^{'})(\alpha,\beta)\]
is a cycle. Finally, if $\alpha\in Y\cup T$, a cycle in $(GoH)[V^{*}\cup \{(\alpha,\beta)\}]$ exists by a reason similar to what was said in Theorem 3.2.

Let each maximal forest of G be without any isolated vertex, and $H$ has a maximal independent set $S$ of size 1. For any maximal forest $F$ in $G$,
\[|S|(K_{2}(F)+L(F))+K_{2}(F)+L^{'}(F)=f(GoH)\]
and therefore we have $|V(F)|=f(GoH)$. Thus $G$ is well-f-covered and $f(G)=f(GoH)$.

Now, let $G$ contains a maximal forest $F^{'}$ with at least one leaf. For any two maximal independent sets $M_{1H}$ and $M_{2H}$ in $H$, we have:  
\[f(H)I(F^{'})+|M_{1H}|(K_{2}(F^{'})+L(F^{'}))+K_{2}(F^{'})+L^{'}(F^{'})=f(GoH),\]
\[f(H)I(F^{'})+|M_{2H}|(K_{2}(F^{'})+L(F^{'}))+K_{2}(F^{'})+L^{'}(F^{'})=f(GoH).\]
Since $K_{2}(F^{'})+L(F^{'})\neq 0$, $|M_{1H}|=|M_{2H}|$ and therefore $H$ is well-covered.
\end{proof}
\end{thm}

\begin{example}
As an application of Theorem 3.5, we shoe that for the graph $G$ with representation:
\begin{center}
\definecolor{ududff}{rgb}{0.30196078431372547,0.30196078431372547,1.}
\begin{tikzpicture}[line cap=round,line join=round,>=triangle 45,x=1.0cm,y=1.0cm]
\clip(7.,0.8) rectangle (9.8,6.3);
\draw [line width=1.6pt] (8.5,5.06)-- (7.7,4.);
\draw [line width=1.6pt] (7.7,4.)-- (9.32,3.98);
\draw [line width=1.6pt] (8.5,5.06)-- (9.32,3.98);
\draw [line width=1.6pt] (9.32,3.98)-- (9.34,2.54);
\draw [line width=1.6pt] (7.7,4.)-- (7.72,2.56);
\draw [line width=1.6pt] (7.72,2.56)-- (9.34,2.54);
\draw (8.26,5.48) node[anchor=north west] {$e$};
\draw (7.24,4.35) node[anchor=north west] {$a$};
\draw (7.28,2.74) node[anchor=north west] {$b$};
\draw (9.35,2.67) node[anchor=north west] {$c$};
\draw (9.34,4.38) node[anchor=north west] {$d$};
\draw (7.8,2) node[anchor=north west] {\small{Figure 1}};
\begin{scriptsize}
\draw [fill=ududff] (8.5,5.06) circle (1.5pt);
\draw [fill=ududff] (7.7,4.) circle (1.5pt);
\draw [fill=ududff] (9.32,3.98) circle (1.5pt);
\draw [fill=ududff] (9.34,2.54) circle (1.5pt);
\draw [fill=ududff] (7.72,2.56) circle (1.5pt);
\end{scriptsize}
\end{tikzpicture}
  \end{center}
  \label{Figure 1}
the composition graph $GoC_{4}$ is not well-f-covered.

The graph $G$ is well-covered and $G[\{a,b,c\}]$ is a maximal forest in $G$, with at least one leaf. The graph $C_{4}$ is well-f-covered and well-covered and has not any maximal independent set of size 1. Also,
\[f(GoC_{4})=6=\alpha(G)f(C_{4}).\]
Thus the conditions (1), (2), and (3)  in Theorem 3.5 hold. Now, consider two maximal forests:
\[F_{1}:=G[\{a,b,c\}], \ \ \ \ F_{2}:=G[\{e,b,c,d\}]\]
in the graph $G$. We have:
\[f(H)I(F_{1})+\alpha(H)(K_{2}(F_{1})+L(F_{1}))+K_{2}(F_{1})+L^{'}(F_{1})=5,\]
\[f(H)I(F_{2})+\alpha(H)(K_{2}(F_{2})+L(F_{2}))+K_{2}(F_{2})+L^{'}(F_{2})=6\]
and therefore the condition (4) is not established; Hence by Theorem 3.5, $GoC_{4}$ is not well-f-covered
\end{example}

The following example shows that the inverse of Theorem 3.5, is not established.

\begin{example}
Consider the graphs $C_{5}$ and $C_{4}$. The graph $C_{5}$ is well-covered and contains a maximal forest with at least one leaf. $C_{4}$ is well-f-covered and well-covered. Also,
\[f(C_{5}oC_{4})=6=\alpha (C_{5})f(C_{4}).\]
Let $F$ be a maximal forest of $C_{5}$. Since $F$ is an induces $P_{4}$, we have:
\[f(C_{4})I(F)+\alpha(C_{4})(K_{2}(F)+L(F))+K_{2}(F)+L^{'}(F)=6=f(C_{5}oC_{4}).\]
Therefore the conditions (1), (2), (3), and (4) in Theorem 3.5 hold. But $C_{5}oC_{4}$ is not well-f-covered because considering the representations:
\begin{center}
\definecolor{ududff}{rgb}{0.30196078431372547,0.30196078431372547,1.}
\begin{tikzpicture}[line cap=round,line join=round,>=triangle 45,x=1.0cm,y=1.0cm]
\clip(3.,1.6) rectangle (11.,5.8);
\draw [line width=1.2pt] (3.84,4.64)-- (5.,4.76);
\draw [line width=1.2pt] (5.,4.76)-- (5.42,3.74);
\draw [line width=1.2pt] (5.42,3.74)-- (4.6,3.);
\draw [line width=1.2pt] (4.6,3.)-- (3.62,3.54);
\draw [line width=1.2pt] (3.84,4.64)-- (3.62,3.54);
\draw [line width=1.2pt] (10.02,4.56)-- (8.66,4.56);
\draw [line width=1.2pt] (8.66,4.56)-- (8.66,3.3);
\draw [line width=1.2pt] (8.66,3.3)-- (10.02,3.3);
\draw [line width=1.2pt] (10.02,4.56)-- (10.02,3.3);
\draw (3.3,5.1) node[anchor=north west] {$x_{1}$};
\draw (5.02,5.12) node[anchor=north west] {$x_{2}$};
\draw (5.47,3.95) node[anchor=north west] {$x_{3}$};
\draw (2.97,3.67) node[anchor=north west] {$x_{5}$};
\draw (4.28,2.98) node[anchor=north west] {$x_{4}$};
\draw (8.11,5.0) node[anchor=north west] {$y_{1}$};
\draw (10.04,4.98) node[anchor=north west] {$y_{2}$};
\draw (10.0,3.36) node[anchor=north west] {$y_{3}$};
\draw (8.11,3.35) node[anchor=north west] {$                y_{4}$};
\draw (9.07,2.66) node[anchor=north west] {$C_{4}$};
\draw (4.24,2.58) node[anchor=north west] {$C_{5}$};
\begin{scriptsize}
\draw [fill=ududff] (3.84,4.64) circle (1.5pt);
\draw [fill=ududff] (5.,4.76) circle (1.5pt);
\draw [fill=ududff] (5.42,3.74) circle (1.5pt);
\draw [fill=ududff] (4.6,3.) circle (1.5pt);
\draw [fill=ududff] (3.62,3.54) circle (1.5pt);
\draw [fill=ududff] (10.02,4.56) circle (1.5pt);
\draw [fill=ududff] (8.66,4.56) circle (1.5pt);
\draw [fill=ududff] (8.66,3.3) circle (1.5pt);
\draw [fill=ududff] (10.02,3.3) circle (1.5pt);
\draw (6.26,2) node[anchor=north west] {\small{Figure 2}};
\end{scriptsize}
\end{tikzpicture}
  \end{center}
  \label{Figure 2}
\vskip 0.6 true cm
the subgraphs $(C_{5}oC_{4})[\{(x_{1},y_{1}),(x_{1},y_{3}),(x_{2},y_{1}),(x_{3},y_{1}),(x_{4},y_{1}),(x_{4},y_{3})\}]$ and $(C_{5}oC_{4})[\{(x_{1},y_{1}),(x_{1},y_{3}),(x_{2},y_{1}),(x_{3},y_{1}),(x_{3},y_{2})\}]$ are two maximal forests of $C_{5}oC_{4}$, of orders 6 and 5, respectively.
\end{example}
 
\vskip 0.4 true cm


\vskip 0.4 true cm




\end{document}